\documentclass[12pt]{article}
\usepackage{amsmath,amsfonts}

\usepackage{graphicx}

\title{Statistical analysis of 
stochastic resonance with ergodic diffusion noise}

\author{Stefano Maria Iacus\\ Department of Economics\\ University of 
Milan\\ Via Conservatorio 7, I-20123 Milan - Italy.\\
E-mail : stefano.iacus@unimi.it
}

\date{November 2001}

\def\dR{\mathbb R}
\def\bE{\mathbf E}
\def\bP{\mathbf P}
\def\de{{\rm d}}
\def\ve{\varepsilon}
\def\vf{\varphi}
\def \mN{{\mathcal N}}
\def \nh{{\mathbf H_0}}  
\def \ah{\mathbf H_1}  
\def\OU{Ornstein-Uhlenbeck}

\begin{document}
\maketitle

\begin{abstract}
A subthreshold signal is transmitted through a channel and may be 
detected when some noise -- with 
known structure and proportional to some level -- is added to the 
data. 
There is an optimal noise level, called 
stochastic resonance, that corresponds to the highest Fisher 
information in the problem of estimation of the signal. As noise we 
consider an ergodic 
diffusion process and the asymptotic is considered as time goes to 
infinity.

We propose consistent
estimators of the subthreshold signal and we solve further a problem 
of hypotheses testing. We also discuss evidence of stochastic 
resonance 
for both estimation and hypotheses testing problems via examples.

\medskip

\noindent
{\bf keywords:} stochastic resonance, diffusion processes, unobservable signal detection, maximum a posterior probability.


\end{abstract}


\section{Introduction}
The term `stochastic resonance' was introduced in the early `80s 
(see \cite{bsv} and \cite{nic}) in the 
study of periodic advance of glaciers on Earth. The 
stochastic resonance is the effect of nonmonotone dependence of the 
response of a system on the noise when this noise  (for instance the 
temperature) is added to a 
periodic input signal (see e.g. \cite{klim}, in which the 
author explains also differences and similarities with the 
notion of stochastic filtering). An extensive review on stochastic 
resonance 
and its presence in different fields of applications can be found in 
\cite{ghjm}.
Following \cite{cbrv}, as stochastic 
resonance we intend  the phenomenon in which the 
transmission of a signal can be improved 
(in terms of statistical quantities)
by the addition of noise. From the statistical point of view the 
problem is  to estimate a signal $\{f(t), \,t<T\}$ transmitted 
through a channel. 
This signal has to be detected by a receiver that can reveal 
signals louder than a threshold $\tau$. If $f(\cdot)$ is bounded from 
above by $\tau$, the signal is not observable and the problem has not 
a solution. But, if some noise $\{\ve(t), \,t<T\}$ is added to the 
signal, the perturbed signal $y(t) = f(t) + \ve(t)$ may be observable 
and inference can be done on $f(\cdot)$. Too few noise is not 
sufficient to give good estimates and too much noise deteriorates 
excessively the signal. The optimal -- in some sense -- level of the 
noise will be called stochastic resonance in this framework. 
Usually (see \cite{ghjm}) the criterion applied to 
measure 
optimality of estimators
are the Shannon mutual information or the Kullback divergence. More 
recently the Fisher information quantity have been also  proposed 
(see \cite{cbrv} and \cite{gww}). Here we are 
concerned  with the Fisher information quantity. It happens that this 
quantity, as a function of the noise, can be maximized for certain 
noise structures. 
If there is only one global maximum, the corresponding noise level is 
the value for which we have stochastic resonance, if several local 
maxima are present, the phenomenon is called stochastic 
multi-resonance.

In this paper we study the problem of estimation and hypotheses 
testing for the following model: we suppose to have a
threshold $\tau>0$ and a subthreshold constant and non negative 
signal 
$\theta$, $\theta<\tau$. We add, in continuous time, a noise 
that is a trajectory of a diffusion process $\{X_t,\, t<T\}$ and we 
observe the perturbed signal $\{Y_t=\theta + \ve X_t,\, t<T\}$ where 
$\ve>0$ is the level of the noise. We propose two schemes of 
observations: {\it i)}  we observe only the proportion of time spent 
by 
the perturbed signal over the threshold $\tau$ and {\it ii)} we 
 measure the energy of the perturbed signal when it is above the 
threshold. The asymptotic is considered as time goes to infinity.

This approach differs from the ones in the current statistical 
literature basically  for two reasons: 
the noise structure is an ergodic diffusion process and not a 
sequence of 
independent and identically distributed random variables and data are 
collected in continuous time. This second 
aspect is a substantial difference but it is not a problem from the
point of view of applications for the two schemes of observations 
proposed if one thinks at analogical devices. 

We propose two different estimators for the schemes 
and we study their asymptotic properties.
We present an example where, in both cases, it emerges the 
phenomenon of stochastic resonance. For the same model we 
also solve the problem of testing the simple hypothesis 
$\theta=\theta_0$ against the simple alternative $\theta=\theta_1$ by 
applying the bayesian maximum a posterior probability criterion. It 
emerges that the overall probability of error is nonmonotonically 
dependent on $\ve$. We show again  that there exists a 
non trivial local minimum of this probability that is again the effect 
of stochastic resonance. The presence of stochastic resonance in 
this context is  noted for the first time here.

The paper is organized as follows.
In Section \ref{sec:model} we set up the regularity 
assumptions of the model. In sections \ref{sec:time} and 
\ref{sec:energy}
we prove some asymptotic properties  estimators for 
the two schemes and  we calculate numerically the points where the 
Fisher information quantity attains its maximum for both models. It 
turns 
out that the 
estimators proposed are asymptotically equivalent to the maximum 
likelihood estimators. Section 
\ref{sec:test} is devoted to the problem of hypotheses testing.  All the figures are collected at the end of the paper.

\section{The model and the structure of the noise}\label{sec:model}
Let $\tau$ be the threshold and $\theta$ a constant signal. Taking 
$0\leq\theta<\tau$ will not influence the calculations that follows 
but may improve the exposition, so we use this assumption. Let 
$\{X_t,t<T\}$ be a given diffusion process solution to the following 
stochastic differential equation
\begin{equation}
\de X_t = S(X_t)\de t + \sigma(X_t) \de W_t
\label{eq:sde}
\end{equation}    
with non random initial value $X_0=0$. The process 
$\{X_t,t<T\}$ is supposed to have the ergodic property 
with  
invariant measure $\bP^*$ and invariant distribution function 
$F(x)=\bP^*((-\infty,x])$ as $T\to\infty$. The functions $S(\cdot)$ 
and $\sigma(\cdot)$ satisfy the global Lipschitz condition
\begin{equation}
|S(x)-S(y)|+|\sigma(x)-\sigma(y)|< L |x-y|,
\tag{C1}
\label{cond:c1}
\end{equation}    
where $L$ is the Lipschitz constant.
Under condition \ref{cond:c1}, equation \eqref{eq:sde} has a unique 
strong solution (see e.g. \cite{ls}) but any equivalent condition to \ref{cond:c1} can 
be assumed because we do not use explicitly it in the sequel.
The following conditions are needed to ensure the ergodicity of the 
process $\{X_t,t<T\}$. If
\begin{equation}
\lim_{|y|\to\infty} \int\limits_0^y \frac{S(u)}{\sigma(u)^2}\de u = 
-\infty
\tag{C2}
\label{cond:c2}
\end{equation}    
and
\begin{equation}
G = \int\limits_\dR \sigma(y)^{-2} \exp\left\{2\int\limits_0^y 
\frac{S(u)}{\sigma(u)^2} \de u\right\}\de y < \infty
\tag{C3}
\label{cond:c3}
\end{equation}    
then there exists the stationary 
distribution function $F(\cdot)$ and it takes the following form
$$
F(x) = G^{-1} \int\limits_{-\infty}^x \sigma(y)^{-2} 
\exp\left\{2\int\limits_0^y 
\frac{S(u)}{\sigma(u)^2} \de u\right\}\de y\,.
$$
Again, any other couple of conditions that imply the existence of $F(\cdot)$ can be used instead of  \ref{cond:c3} and \ref{cond:c3}.

We perturb the signal $\theta$ by adding, proportionally to some 
level $\ve>0$, the trajectory diffusion process $\{X_t,t<T\}$ into the channel. 
The result will be the perturbed 
signal $\{Y_t = \theta + \ve\,X_t,t<T\}$. This new signal will be 
detectable only when it is above the threshold $\tau$. Moreover, 
$\{Y_t,t<T\}$
is still ergodic with trend and diffusion coefficients respectively 
$S_\theta(y) = S((y-\theta)/\ve)$ and $\sigma_\theta((y-\theta)/\ve)$ 
and initial value $Y_0 = 
\theta$, but we will not use directly this process. We denote by
$\{M_t = Y_t\,\chi_{\{Y_t>\tau\}},t<T\}$ the observable part of the 
trajectory of $\{Y_t,t<T\}$, being $\chi_A$ the indicator function of 
the set $A$.

We consider two possibile schemes of observation:
\begin{itemize}
\item[{\it i)}] we observe only the proportion of time spent by 
$\{Y_t,t<T\}$  over the threshold $\tau$
$$
\Gamma_T = \frac1{T} \int\limits_0^T \chi_{\{Y_t>\tau\}} \de t\,,
$$
\item[{\it ii)}] we measure the energy of the signal 
$\{M_t,t<T\}$  
$$
\nu_T = \frac1{T} \int\limits_0^T M_t^2 \de t\,.
$$
\end{itemize}
In the next sections, for the two models we establish asymptotic 
properties of estimators given by the generalized method of moments.
In \cite{k2000} different properties of the generalized method of 
moments for ergodic diffusion processes are studied. In this note we 
follows the lines given in the paper of\cite{gww} for the i.i.d. setting.
These results are interesting in themselves independently from the 
problem of stochastic resonance.
We give an example of stochastic resonance based on the \OU\, 
process where the phenomenon of stochastic resonance appears 
pronounced and in which results in a closed form can be written down.

\section{Observing the time spent by the process over the 
threshold}\label{sec:time}
The random variable $\Gamma_T$ can be rewritten in terms of the 
process $\{X_t,t<T\}$ as
$$
\Gamma_T = \frac{1}{T} \int_0^T 
\chi_{\left\{X_t>\frac{\tau-\theta}{\ve}\right\}} \de t\,.
$$
By the ergodic property of $\{X_t,t<T\}$ we have that
\begin{equation}
\Gamma_T \underset{T\to\infty}{\longrightarrow} \bE 
\chi_{\left\{\xi>\frac{\tau-\theta}{\ve}\right\}}
=\bP\left(\xi>\frac{\tau-\theta}{\ve}\right) = 1 - 
F\left(\frac{\tau-\theta}{\ve}\right)
=\pi
\label{eq:teta}
\end{equation}
where $\xi$ has $F(\cdot)$ as distribution function.
From \eqref{eq:teta} it derives that
$$
\theta = \theta(\pi) = \tau + \ve F^{-1}(1-\pi)\,,
$$
so that $\theta$ is a one-to-one continuous function of $\pi$.
From the Glivenko-Cantelli theorem (see e.g. \cite{k2000}) for the 
empirical distribution 
function (EDF) defined by 
$$
\hat F_T(x) = \int_{-\infty}^x \chi_{\{X_t<x\}}\de t
$$
follows directly that $\Gamma_T$ is a $\sqrt 
T$-consistent estimator of $\pi$ thus also
$$
\hat \theta_T = \theta(\Gamma_T) = \tau + \ve F^{-1}(1-\Gamma_T)\,,
$$
is a $\sqrt T$-consistent estimator for $\theta$. We can calculate 
the asymptotic variance of this estimator. It is known that (see \cite{k1997a} and \cite{ilia}) the
 EDF is asymptotically Gaussian and in particular
$$
\sqrt{T}\left(\hat F_T(x)-F(x)\right)\Longrightarrow \mN\left(0, 
V(x)\right)\,,
$$
where $V(x)=I_F(x)^{-1}$
 is the inverse of the analogue of the Fisher information 
quantity  in the problem of 
distribution function estimation :
\begin{equation}
I_F(x) = \left(
4 \bE \left(
\frac{F(\xi\wedge x)(1-F(\xi\vee x))}{\sigma(\xi)f(\xi)}
\right)^2
\right)^{-1}\,,
\label{eq:fisher}
\end{equation}
where $a \wedge b = \min(a,b)$ and $a \vee b=\max(a,b)$.
The quantity $V(x)$ is also the minimax asymptotic lower bound for 
the 
quadratic risk associated to the estimation of $F(x)$, so that $\hat 
F_T(x)$ 
is asymptotically efficient in this sense.

The asymptotic variance $\Sigma(\theta)$  of $\hat\theta_T$ can be 
derived by means of 
the so-called $\delta$-method (see e.g. \cite{bncox}):
$$
\sqrt{T}\left(\theta\left(\Gamma_T\right)-\theta(\pi)\right) =
\sqrt{T}\left(\Gamma_T-\pi\right)\theta^\prime(\pi)+o_T\left(|\Gamma_T-\pi|\right)
$$
thus
$$
\Sigma(\theta) = 
\ve^2\frac{V\left(\frac{\tau-\theta}{\ve}\right)}{f\left(\frac{\tau-\theta}{\ve}\right)^2}\,,
$$
where $f(\cdot)$ is the density of $F(\cdot)$.
The quantity $\Sigma(\theta)$ can also be derived from the asymptotic 
minimal variance $V(\cdot)$ 
of the EDF estimator. In 
fact, with little abuse of notation,
by putting $F(x) = \mu$ and $\theta(\mu) = \tau - \ve F^{-1}(\mu)$ we 
have that
$$
V(\theta(\mu)) = \frac{(\theta'(\mu))^2}{I_F(\theta(\mu))} = 
\Sigma(\theta)\,.
$$
\subsection{Link with the likelihood estimator}
We now show that $\Gamma_T$  also maximizes the approximate 
likelihood of the model. In fact, for the central limit theorem for 
the EDF we have
$$
\hat F_T(x) = F(x) + \sqrt{\frac{V(x)}{T}} Z + o_T(1)
$$
where $Z$ is a standard Gaussian random variable. Thus
$$
\begin{aligned}
\bP\left(
\Gamma_T< \gamma
\right) 
&= \bP\left(
1 - \hat F_T\left(\frac{\tau-\theta}{\ve}\right)
<\gamma
\right)\\
&=
1 - \Phi\left(  
\left(1-\gamma-F\left(\frac{\tau-\theta}{\ve}\right)\right) 
\sqrt{\frac{T}{V\left(\frac{\tau-\theta}{\ve}\right)}}
\right)+o_T(1)
\end{aligned}
$$
where $\Phi()$ is the distribution function of $Z$.
We approximate the likelihood function of $\Gamma_T$ by  
$$
\vf(\theta;\Gamma_T) \simeq  \frac{\sqrt{T}}{\sqrt{2\pi 
V\left(\frac{\tau-\theta}{\ve}\right)}}
\exp\left\{
-\frac{T}2 
\frac{\left(1-\Gamma_T-F\left(\frac{\tau-\theta}{\ve}\right)\right)^2}{V\left(\frac{\tau-\theta}{\ve}\right)}
\right\}
$$
that is maximal when 
$$
1-\Gamma_T-F\left(\frac{\tau-\theta}{\ve}\right) = 0
$$
thus, the maximum likelihood estimator of $\theta$ (constructed on 
the approximated likelihood) reads
$$
\hat\theta_T = \tau-\ve F^{-1}\left(1-\Gamma_T\right)\,.
$$
So if  the approximation above is acceptable, one can infer the 
optimality property of $\Gamma_T$ of having minimum variance from 
being also the maximum likelihood estimator.
\subsection{An example 
of stochastic resonance}
To view the effect of stochastic resonance on the Fisher information 
we consider a particular 
example. By setting $S(x) = -x$ and $\sigma(x) = 1$ the noise 
become a standard \OU\, process solution to the stochastic 
differential equation
$$
\de X_t = -X_t \de t + \de W_t\,.
$$
In such a case, the ergodic distribution function $F(\cdot)$ is the 
Gaussian 
law with zero mean and variance 1/2.
The asymptotic variance $\Sigma(\theta)$ assumes the following 
form
$$
\Sigma(\theta) = 
\frac{\ve^2\, 
\pi^\frac32}{e^{-2\left(\frac{\tau-\theta}{\ve}\right)^2}} 
\int_{\dR} 
\left(1+{\rm erf}\left(x\wedge 
\frac{\tau-\theta}{\ve}\right)\right)^2\left(
1-{\rm erf}\left(x\vee 
\frac{\tau-\theta}{\ve}\right)
\right)^2
e^{x^2}
\de x
$$
where ${\rm erf}(x) = \frac2\pi \int_0^x e^{-t^2}\de t$ is the 
classical error function.

\begin{figure}[h]
\begin{center}
\begin{tabular}{c }
\includegraphics{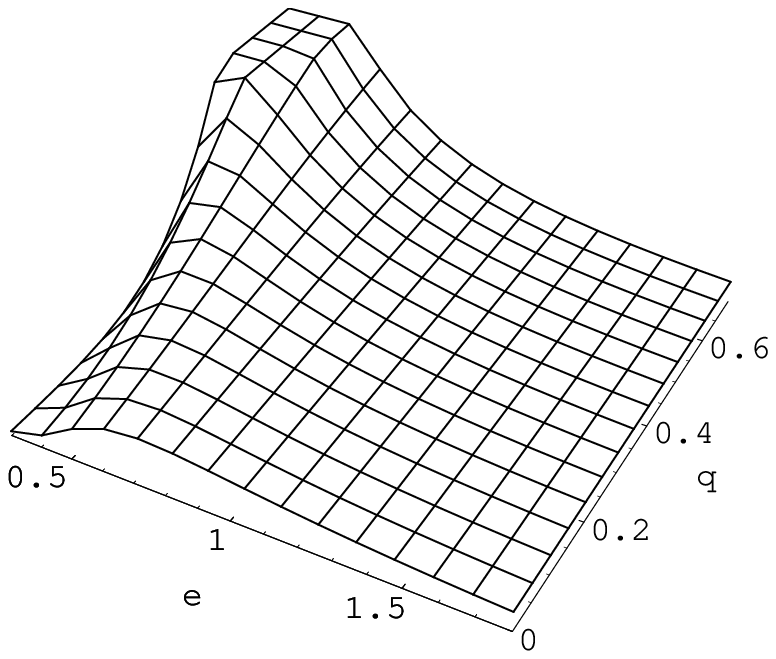}   \\ 
\includegraphics{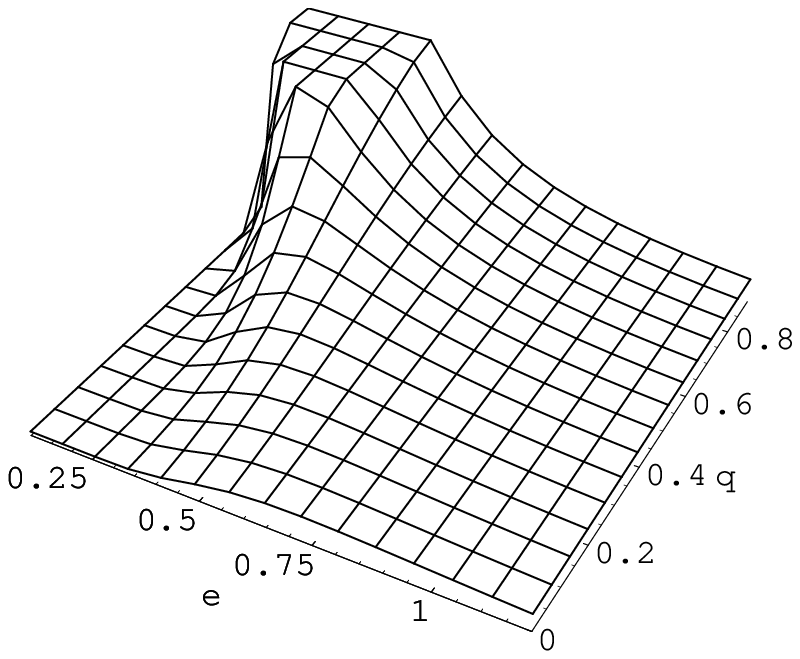}    \\
\end{tabular}
\end{center}
\caption{The Fisher information $\Sigma(\theta)^{^-1}$ when 
observing $\Gamma_T$ (up) and 
$\tilde\Sigma(\theta)^{^-1}$ when the observation is $\nu_T$ (down) 
exhibit 
stochastic resonance. The threshold is fixed at $\tau=1$.}
\label{fig1}
\end{figure}

In Figure \ref{fig1} it is shown that for this model there exists the 
phenomenon of stochastic resonance. For a fixed level of noise $\ve$ 
the Fisher information increases as the signal $\theta$ is closer to 
the 
threshold $\tau$. For a fixed value of the signal $\theta$, the 
Fisher 
information, as a function of $\ve$, has a single maximum, that is the
optimal level of noise. For example, if $\theta=0$ then the optimal 
level is $\ve^*=0.1811$ and for $\theta = 0.5$ it is $\ve^*=0.7244$.

\section{Measuring the observed energy}\label{sec:energy}
Suppose that it is possibile to observe not only the time when the 
perturbed process is over the threshold but also its trajectory above 
$\tau$, say $M_t = Y_t \chi_{\{Y_t>\tau\}}$, $t<T$.
We now show how it is possibile to estimate the unknown signal 
$\theta$ from the equivalent of the energy of the signal for  $M_t$: 
literally from the quantity
$$
\nu_T=\frac1{T} \int_0^T M_t^2 \de t\,.
$$
We use the following general  result from \cite{k1997b} on the 
estimation 
of functionals of the invariant distribution functions for ergodic 
diffusion processes.
{\it Let $R(\cdot)$ and $N(\cdot)$ be such that
$\bE \left(|R(\xi)S(\xi)|+|N(\xi)|\right)<\infty$. Then
$$
\frac1{T} \int_0^T R(X_t) \de X_t + \frac1{T} \int_0^T N(X_t) \de 
t$$
is a $\sqrt T$-consistent estimator for 
$\nu=\bE\left(R(\xi)S(\xi)+N(\xi)\right)$
where $\xi$ is distributed according to $F(\cdot)$}. 
In our case $R(\cdot)=0$ and $N(x)=(\ve x 
+\theta)^2\chi_{\{x>\frac{\tau-\theta}{\ve}\}}$.
The estimator $\nu_T$ can be rewritten as
$$
\nu_T=\frac1{T}
\int_0^T \left(\ve\,X_t + \theta\right)^2 
\chi_{\{x>\frac{\tau-\theta}{\ve}\}}\de t
$$
and it converges to the quantity
$$
\nu=\nu(\theta) = \ve^2\,\bE\left(\xi^2 
\chi_{\{x>\frac{\tau-\theta}{\ve}\}}\right)+
\theta^2\left(1-F\left(\frac{\tau-\theta}{\ve}\right)\right)
+2\,\theta\,\ve \,\bE\left(\xi \,
\chi_{\{x>\frac{\tau-\theta}{\ve}\}}\right)
$$
that is a continuous and  increasing function of $\theta$, 
$0<\theta<\tau$. 
Its inverse $\theta(\nu)=\nu^{-1}(\nu)$ allows us to have again
$\tilde\theta_T = \theta(\nu_T) = \nu^{-1}(\nu_T)$.
By applying the 
$\delta$-method once again, we can obtain the asymptotic 
variance of $\tilde\theta_T$ from the asymptotic variance of $\nu_T$.
$$
\sqrt{T}\left(\theta(\nu_T)-\theta(\nu)\right) = 
\frac{\sqrt{T}\left(\nu_T-\nu 
\right)}{\nu'(\theta)}+o_T\left(|\nu_T-\nu|\right)
$$
where
$$\nu'(\theta)=
\frac{\tau^2}{\ve}\, f\left(\frac{\tau-\theta}{\ve}\right) + 2  
\, \theta \, \left(1-F\left(\frac{\tau-\theta}{\ve}\right)\right)+2\, 
\ve\,\bE\left(\xi \, 
\chi_{\left\{x>\frac{\tau-\theta}{\ve}\right\}}\right)\,.
$$
The asymptotic variance of $\nu_T$ is given by (see \cite{k1997b}) 
$$\tilde V(\theta) = \tilde V(\nu(\theta) ) = 4\,\bE\left\{ 
\frac{M(\xi)^2}{f(\xi)^2}\right\}\,,$$ 
where
$$M(y) = \bE\left\{\left(F(y)-\chi_{\{\xi<y\}}\right)
(\ve \xi +\theta)^2\chi_{\{\xi>\frac{\tau-\theta}{\ve}\}}
\right\}
$$
and its inverse is also 
the minimal asymptotic variance in the problem of estimation of 
functionals for ergodic diffusion.
Thus, the asymptotic variance of $\tilde\theta_T$ is given by
$$\tilde\Sigma(\theta) = 4\,\left.\bE\left\{ 
\frac{M(\xi)^2}{f(\xi)^2}\right\}\right/\nu'(\theta)^2\,.$$

{\bf Remark:} \, {\sl
By the asymptotic normality of $\nu_T$ follows that $\tilde\theta_T$ 
is also the  value that maximizes the approximate likelihood 
function. In 
fact, as in the previous example, if we approximate the density 
function  of $\nu_T$ with
$$
\vf(\theta;\nu_T) \simeq \frac{\sqrt{T}}{\sqrt{2\pi \tilde 
V\left(\frac{\tau-\theta}{\ve}\right)}}
\exp\left\{
-\frac{T}2 
\frac{\left(\nu_T-\nu(\theta)\right)^2}{\tilde 
V\left(\frac{\tau-\theta}{\ve}\right)}
\right\}
$$
it is clear that  $\tilde\theta_T$ is its maximum.
}

\subsection{The effect of stochastic resonance}
As before, we put in evidence the phenomenon of stochastic resonance 
by using the  \OU\, process as noise.
The quantities involved ($\nu(\theta)$ and $\nu'(\theta)$) transform 
into the following
$$
\nu(\theta)=\frac14 \, \Bigg\{\, \, \ve^2+2\, \theta^2+
2\, \ve \, 
\frac{\theta +\tau}{\sqrt{\pi}} 
e^{-\left(\frac{\tau-\theta}{\ve}\right)^2}
-\, \, (\ve^2+2\, \theta^2)\,
{\rm erf}\left(\frac{\tau-\theta}{\ve}\right)\Bigg\}
$$
(from which it appears that $\nu(\theta)$ is an increasing function 
of 
$\theta$) and
$$
\nu'(\theta) = \theta 
+\frac{(\ve^2+\tau^2)\,e^{-\left(\frac{\tau-\theta}{\ve}\right)^2}}{\ve\,\sqrt{\pi}}
-\theta\,{\rm erf}\left(\frac{\tau-\theta}{\ve}\right)\,.
$$
In Figure \ref{fig1} it is plotted the Fisher information of the 
model as a function of $\theta$ and $\ve$. Also in this case there is 
evidence
of stochastic resonance. For a fixed value of 
$\theta$ is then possibile to find the optimal noise level $\ve$. For 
example, taking $\theta = 0$ then we have stochastic resonance at 
$\ve^* = 0.7234$ and for $\theta=0.5$, $\ve^* =0.3636$.

\section{Hypotheses testing problem}\label{sec:test}
We now study a problem of testing two simple hypotheses for the 
model discussed in the previous section. As in \cite{cb}, we apply the maximum {\it a posteriori} probability (MAP) 
criterion. We will see that the decision rules for our model are
similar to the one proposed by Chapeau-Blondeau in the i.i.d. setting.

Given the observation $\Gamma_T$ we want to verify the null 
hypothesis that 
the unknown constant signal is $\theta_0$ against the simple 
alternative $\theta_1$, with $\theta_0<\theta_1<\tau$:
$$
\begin{aligned}
\nh &: \theta = \theta_0\\
\ah &: \theta = \theta_1
\end{aligned}    
$$
Suppose that, before observing $\Gamma_T$, we have a {\em prior} 
information
on the parameter, that is  $\bP_0 = 
\bP(\theta=\theta_0)$ and $\bP_1 = 
\bP(\theta=\theta_1)$.
The MAP criterion uses the following likelihood ratio
$$
\lambda = 
\frac{\bP(\theta=\theta_1|\Gamma_T)}{\bP(\theta=\theta_0|\Gamma_T)}
=
\frac{\vf(\theta_1;\Gamma_T)\,\bP_1}{\vf(\theta_0;\Gamma_T)\,\bP_0}
$$
and the decision rule is to accept $\nh$ whenever $\lambda>1$ 
(decision $D_1$) or refuse it otherwise (decision $D_0$).  
The overall probability of error is
$$
\bP_{err} = \bP(D_1|\nh) \bP_0 + \bP(D_0|\ah) \bP_1\,.
$$
Let now be
\begin{equation}
\sigma_i^2 = 
\frac{V\left(\frac{\tau-\theta_i}{\ve}\right)}{T}\,,\qquad
\mu_i = 1-F\left(\frac{\tau-\theta_i}{\ve}\right),\qquad i=0,1\,.
\label{eq:parametri}
\end{equation}
Then, the likelihood $\lambda$ appears as
$$
\lambda = \frac{\sigma_0\,\bP_1}{\sigma_1\,\bP_0}
\exp\left\{
-\frac12 
\left(
\frac{(\Gamma_T-\mu_1)^2}{\sigma_1^2}-
\frac{(\Gamma_T-\mu_0)^2}{\sigma_0^2}\right)
\right\}
+o_T(1)
$$
To write explicitly the decision rule and then study the effect of 
stochastic 
resonance we have to distinguish three cases: $\sigma_0>\sigma_1$, 
$\sigma_0<\sigma_1$ and $\sigma_0=\sigma_1$.
\begin{enumerate}
\item    Let it be $\sigma_0>\sigma_1$, then put
$$
\begin{aligned}
\sigma_2 &= \sqrt{\sigma_0^2-\sigma_1^2} &\phantom{cccc} 
\gamma' &= 
\frac{\mu_1\,\sigma_0^2-\mu_0\, 
\sigma_1^2-\sigma_0\,\sigma_1\,\sqrt{\Delta}}{\sigma_2^2} \\
\Delta &= 
(\mu_0-\mu_1)^2 - 2\,\sigma_2^2\, \log \left(
\frac{\bP_0\,\sigma_1}{\bP_1\,\sigma_0}\right) &
\gamma'' &= 
\frac{\mu_1\,\sigma_0^2-\mu_0\, 
\sigma_1^2+\sigma_0\,\sigma_1\,\sqrt{\Delta}}{\sigma_2^2} 
\end{aligned}
$$
Then, if $\Delta<0$ accept $\nh$ and $\bP_{err}=\bP_1$.
If $\Delta>0$, then if $\gamma' < \Gamma_T < \gamma''$ reject $\nh$ 
otherwise accept it. In both cases
$$
\begin{aligned}
\bP_{err} &= \frac12 \left\{
{\rm erf}\left(\frac{\gamma'' - \mu_0}{\sigma_0 \sqrt2}\right)
-{\rm erf}\left(\frac{\gamma' - \mu_0}{\sigma_0 \sqrt2}\right)
\right\}\,\bP_0\\
&\phantom{+} +
\frac12 \left\{2-
{\rm erf}\left(\frac{\gamma'' - \mu_1}{\sigma_1 \sqrt2}\right)
+{\rm erf}\left(\frac{\gamma' - \mu_1}{\sigma_1 \sqrt2}\right)
\right\}\,\bP_1\,.
\end{aligned}
$$
\item     Let it be $\sigma_0<\sigma_1$, then put
$$
\begin{aligned}
\sigma_2 &= \sqrt{\sigma_1^2-\sigma_0^2} &\phantom{cccc} 
\gamma' &= 
\frac{\mu_0\, 
\sigma_1^2-\mu_1\,\sigma_0^2-\sigma_0\,\sigma_1\,\sqrt{\Delta}}{\sigma_2^2} 
\\
\Delta &= 
(\mu_0-\mu_1)^2 - 2\,\sigma_2^2\, \log \left(
\frac{\bP_1\,\sigma_0}{\bP_0\,\sigma_1}\right) &
\gamma'' &= 
\frac{\mu_0\, 
\sigma_1^2-\mu_1\,\sigma_0^2+\sigma_0\,\sigma_1\,\sqrt{\Delta}}{\sigma_2^2} 
\end{aligned}
$$
Then, if $\Delta<0$ reject $\nh$ and $\bP_{err}=\bP_0$.
If $\Delta>0$, then if $\gamma' < \Gamma_T < \gamma''$ accept $\nh$ 
otherwise reject it. In both cases
\begin{equation}
\begin{aligned}
\bP_{err} &= \frac12 \left\{
{\rm erf}\left(\frac{\gamma'' - \mu_1}{\sigma_1 \sqrt2}\right)
-{\rm erf}\left(\frac{\gamma' - \mu_1}{\sigma_1 \sqrt2}\right)
\right\}\,\bP_1\\
&\phantom{+} +
\frac12 \left\{2-
{\rm erf}\left(\frac{\gamma'' - \mu_0}{\sigma_0 \sqrt2}\right)
+{\rm erf}\left(\frac{\gamma' - \mu_0}{\sigma_0 \sqrt2}\right)
\right\}\,\bP_0\,.
\end{aligned}
\label{eq:perr}
\end{equation}
\item    Let it be $\sigma_0=\sigma_1$, then put
$$
\gamma = \frac{\mu_1^2 - \mu_0^2 + 2\,\sigma_0^2\, 
\log\left(\frac{\bP_0}{\bP_1}\right)}{2(\mu_1-\mu_0)}
$$
Then, if $\Gamma_T > \gamma$ reject $\nh$ otherwise 
accept it. In both cases
$$
\bP_{err} = \frac12 \left\{
1+  {\rm erf}\left(\frac{\gamma - \mu_1}{\sigma_1 
\sqrt2}\right)\,\bP_1
-{\rm erf}\left(\frac{\gamma - \mu_0}{\sigma_0 \sqrt2}\right)\,\bP_0
\right\}\,.
$$
\end{enumerate}

\subsection{Example}
\begin{figure}[ht]
\begin{center}
\includegraphics{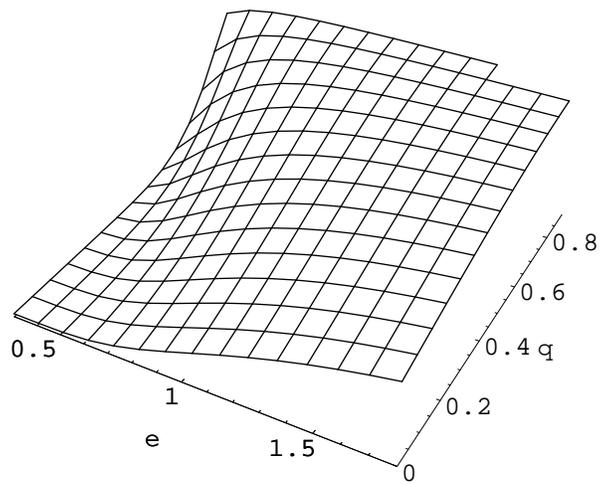}    
\end{center}
\caption{The variance $V((\tau-\theta)/ \ve)$ with $\tau=1$.}
\label{fig2}
\end{figure}

\begin{figure}
\begin{center}
\begin{tabular}{c c c}
\makebox{\includegraphics[width=3cm]{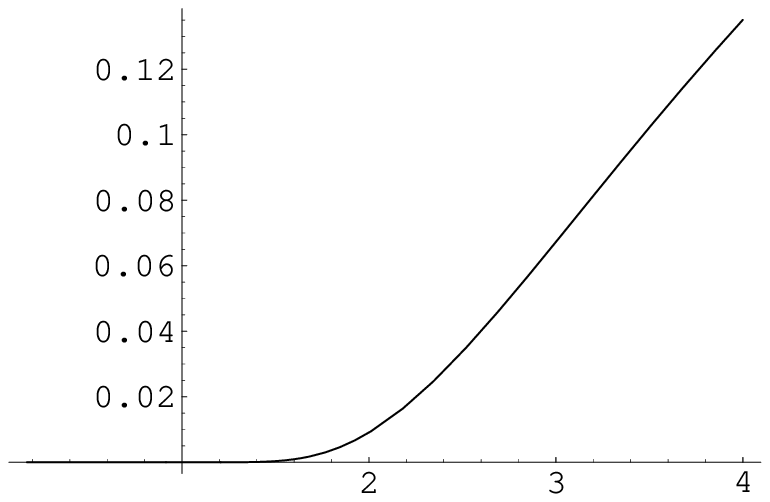}} & 
\includegraphics[width=3cm]{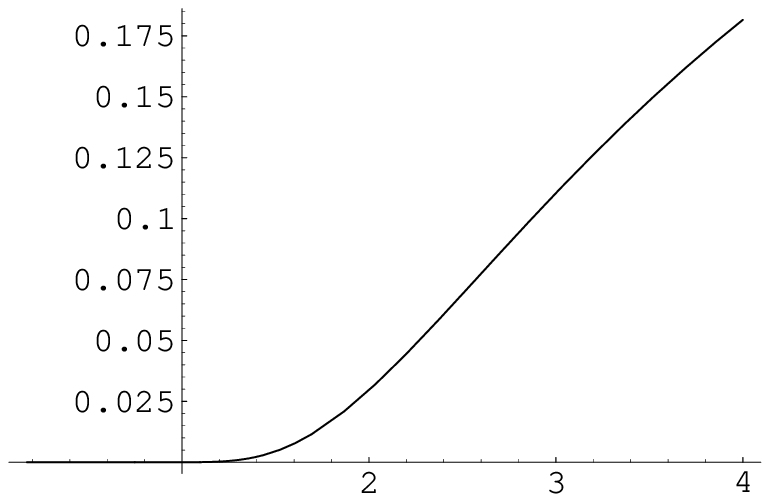} & 
\includegraphics[width=3cm]{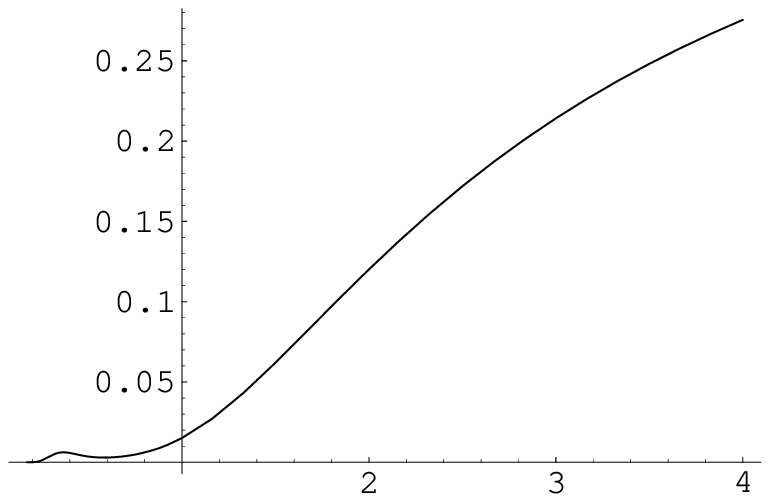} \\   
\includegraphics[width=3cm]{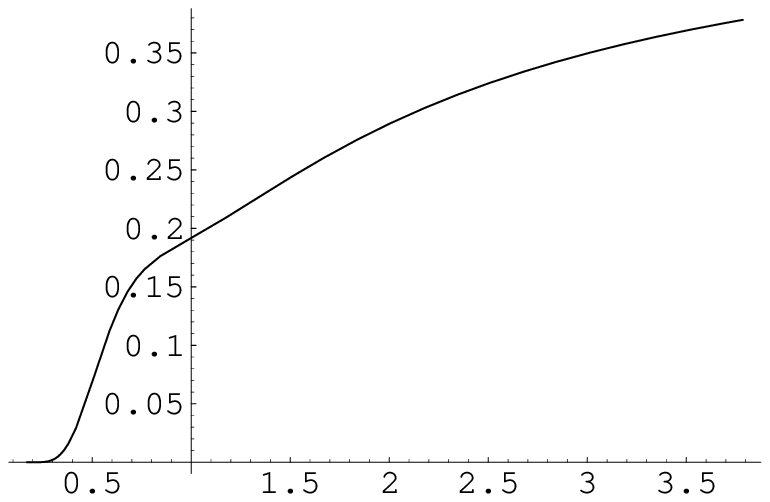} & 
\includegraphics[width=3cm]{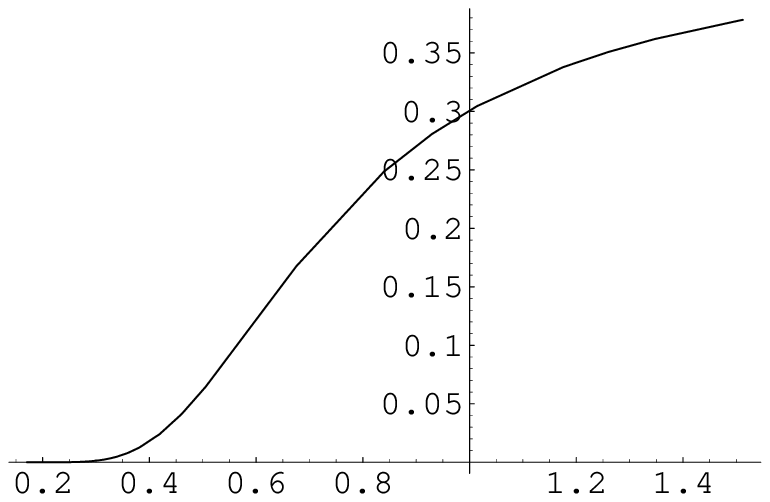} & 
\includegraphics[width=3cm]{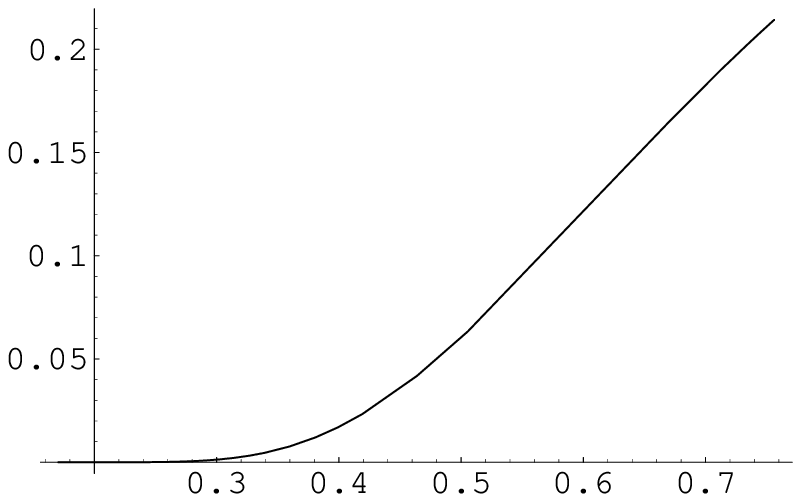}
\end{tabular}    
\includegraphics{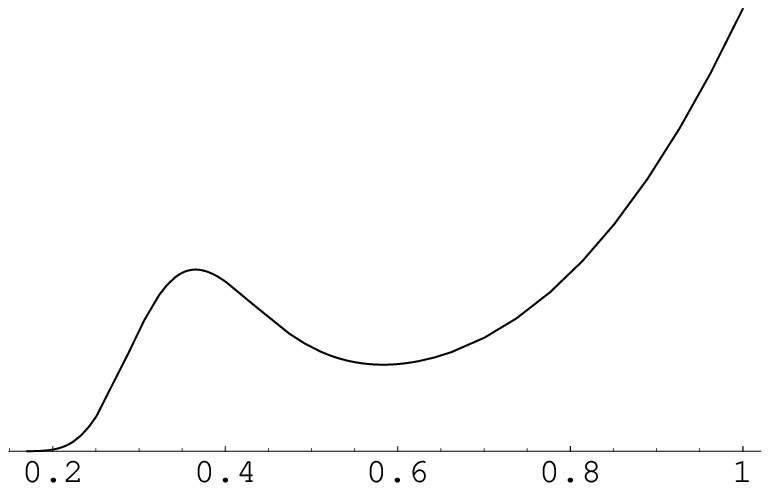}
\end{center}
\caption{The significant part of the plot of $\bP_{err}$ as a 
function of $\ve$ for 
different values of
$\theta_1$ given $\tau=1$ and $\theta_0=0$. From left-top to 
right-bottom: $\theta_1$ = 0.1, 0.25, 0.5, 0.75, 0.9, 0.95. The bottom
plot is to show the effect of stochastic resonance for $\theta_1$ = 
0.5.}
\label{fig3}
\end{figure}

As before, we apply this method to the \OU\, model. In this case the 
variance $\sigma_i^2 = \sigma_i^2(\theta,\ve,\tau,T)$, for a fixed 
threshold $\tau$ and noise level $\ve$, is a non decreasing function 
of $\theta$ being $T$ only a scale factor (Figure \ref{fig2} 
gives a numerical representation of this statement). Thus, the 
$\bP_{err}$ is, in general, given by formula \eqref{eq:perr}.
What is amazing is the behavior of $\bP_{err}$. In Figure \ref{fig3} 
it is reported the graph of $\bP_{err}$ as a function of $\ve$ and
$\theta_1$ given $\tau=1$ and $\theta_0=0$. For $\theta_1$ around 
1/2 the $\bP_{err}$ shows the effect of stochastic resonance. So it 
appears that in some cases the noise level $\ve$ can reduce sensibly 
the overall probability of making the wrong decision. This kind of 
behavior is non outlined in the work of \cite{cb}.

{\bf Remark:} \, {\sl
Following the same scheme, similar results can be obtained for the 
model {\it ii)} 
when we observe the energy $\nu_T$. In this case it is sufficient to 
replace in \eqref{eq:parametri} the values of $\sigma^2_i$ and 
$\mu_i$ with the 
quantities
$$
\sigma^2_i = \frac{\tilde V(\theta)}{T},\qquad \mu_i = 
\nu(\theta_i),\qquad i=0,1
$$
and  $\Gamma_T$ with $\nu_T$ in the decision rule.
}

\section*{Final remarks}
The use of ergodic diffusions as noise in the problem of stochastic 
resonance seems quite powerful. Characterizations of classes of 
ergodic process that enhance the stochastic resonance can be 
done (see e.g. \cite{gww}) but not in a simple way as in the 
i.i.d. case as calculations are always cumbersome. The problem of a 
parametric non constant 
signal can also be treated while the full nonparametric non constant 
signal requires more attention and will be object for further 
investigations. 
For i.i.d. observations, \cite{mu} and \cite{mw} considered the problem of non parametric 
estimation for regression models of the form $Y(t_i) = 
s(t_i)+\sigma(t_i)$, $i=1,\ldots, n$.
Their approach can be applied 
in this context.

Other criterion of optimality than  the Fisher 
information quantity can be used as it is usually done in information 
theory (e.g. Shannon mutual information or Kullback divergence).

The analysis of the overall probability of error seems to put in 
evidence something new with respect to the current literature 
(see e.g. \cite{cb}).
It is worth noting that in a recent paper \cite{bg}  models driven by ergodic diffusions have also been used but 
the effect of stochastic resonance is not used to estimate 
parameters.

\newpage

\section*{Figures}

\end{document}